\documentclass{amsart}

\newcommand{\Z}{\ensuremath{\mathbf Z}}
\newcommand{\N}{\ensuremath{ \mathbf N }}

\newtheorem{theorem}{Theorem}
\newcommand{\bt}{\begin{theorem}}
\newcommand{\et}{\end{theorem}}
\newtheorem{lemma}{Lemma}
\newcommand{\bl}{\begin{lemma}}
\newcommand{\el}{\end{lemma}}
\newcommand{\beq}{\begin{equation}}
\newcommand{\eeq}{\end{equation}}
\newcommand{\benum}{\begin{enumerate}}
\newcommand{\eenum}{\end{enumerate}}
\newcommand{\ord}{\text{ord}}

\title[Sums of products]{Sums of products of congruence classes\\ and of arithmetic progressions}

\author{Sergei V. Konyagin}
\address{Department of Mathematics\\
Moscow State University\\
Moscow, Russia}
\email{konyagin@ok.ru}
\thanks{The work of S.V.K was supported in part by the INTAS grant 03-51-5070.}
\author{Melvyn B. Nathanson}
\address{Department of Mathematics\\
Lehman College (CUNY)\\
Bronx, New York 10468}
\email{melvyn.nathanson@lehman.cuny.edu}
\thanks{The work of M.B.N. was supported in part by grants from the NSA Mathematical Sciences Program and the PSC-CUNY Research Award Program.}

\keywords{Sumsets, product sets, sums of products, congruence classes, arithmetic progressions, linear diophantine problem of Frobenius}
\subjclass[2000]{Primary 11A07, 11B25, 11B75.} 

\date{\today}

\begin{document}

\begin{abstract}

Consider the congruence class $R_m(a)=\{a+im:i\in \Z\}$ and the infinite arithmetic progression $P_m(a)=\{a+im:i\in \N_0\}$.  For positive integers $a,b,c,d,m$ the sum of products set $R_m(a)R_m(b)+R_m(c)R_m(d)$ consists of all integers of the form $(a+im)(b+jm)+(c+km)(d+\ell m)$ for some $i,j,k,\ell\in \Z\}.$ It is proved that if  $\gcd(a,b,c,d,m)=1,$ then $R_m(a)R_m(b)+R_m(c)R_m(d)$ is equal to the congruence class $R_m(ab+cd),$ and that the sum of products set $P_m(a)P_m(b)+P_m(c)P_m(d)$ eventually coincides with the infinite arithmetic progression $P_m(ab+cd).$
\end{abstract}

\maketitle

\section{Sums of product sets}
Let \Z\ denote the set of integers and $\N_0$ the set of nonnegative integers.  For every prime $p$ and integer $n$, we denote by $\ord_p(n)$ the greatest integer $k$ such that $p^k$ divides $n$.

Let $X$ and $Y$ be sets of integers.  These sets \emph{eventually coincide}, denoted $X \sim Y,$ if there is an integer $n_0$ such that, for all $n \geq n_0,$ we have $n \in X$ if and only if $n \in Y.$  We define the \emph{sumset} $X+Y = \{x+y : x\in X, y \in Y \}$, the \emph{product set} $XY = \{xy : x\in X, y \in Y \}$, and, for any integer $\delta,$ the \emph{dilation} $\delta \ast X  = \{\delta x : x\in X \}$.  

Let $a,b,$ and $m$ be integers with $m\geq 1$.  We denote the congruence class of $a$ modulo $m$ by
\[
R_m(a)=\{a+im:i\in \Z\}.
\]
For all $a$ and $b$, we have
\[
R_m(a)R_m(b) \subseteq R_m(ab).
\]
This inclusion can be strict.  For example, $53 \in R_{19}(15)$ but $53 \notin R_{19}(3)R_{19}(5)$ since 53 is prime.  Thus, the product of two congruence classes modulo $m$ is not necessarily a congruence class modulo $m$.

The case of sums of products of congruence classes is different.  For all integers $a,b,c,d,$ and $m$ with $m \geq 1$ we have
\[
R_m(a)R_m(b) + R_m(c)R_m(d) \subseteq R_m(ab+cd).
\]
We shall prove that if $\gcd(a,b,c,d,m)=1$, then   
\[
R_m(a)R_m(b)+R_m(c)R_m(d) = R_m(ab+cd)
\]
and if $\gcd(a,b,c,d,m)=\delta$, then   
\[
R_m(a)R_m(b)+R_m(c)R_m(d) = R_{\delta m}(ab+cd).
\]

Let $a,b,c,d,$ and $m$ be positive integers.  We denote the infinite arithmetic progression with initial term $a$ and difference $m$ by
\[
P_m(a)=\{a+im:i\in \N_0\}.
\]
Then $P_m(a)P_m(b) \subseteq P_m(ab)$.  Again we choose $a=3, b=5,$ and $m=19.$  By Dirichlet's theorem, there are infinitely many primes $p\equiv 15\pmod{19},$ and none of these is a product of an integer congruent to 3 and an integer congruent to 5 modulo 19.  It follows that there are infinitely many integers in the arithmetic progression $P_{19}(15)$ that do not belong to the product set $P_{19}(3)P_{19}(5),$ and so 
 \[
P_{19}(3) P_{19}(5) \not\sim P_{19}(15).
\]
Thus, the product of two arithmetic progressions with difference $m$ does not necessarily eventually coincide with an arithmetic progression with difference $m$.  On the other hand, we shall prove that if $(a,b,c,d,m)=1,$ then 
 \[
P_m(a)P_m(b)+P_m(c)P_m(d) ~\sim P_m(ab+cd).
\]
 
Sums of products in finite fields have been studied recently by Hart and Iosevich~\cite{hart-iose07} and Glibichuk and Konyagin~\cite{glib-kony07}.  The problem of sums of products of congruence classes of integers actually arose in unpublished work of Robert Schneiderman and Peter Teichner in low-dimensional topology.  They are studying the failure of the Whitney move and are trying to measure this failure in terms of an obstruction theory for "Whitney Towers" (iterated layers of Whitney disks) built on immersed surfaces in 4-manifolds. The associated intersection invariants have indeterminacies that can be non-linear in the presence of certain non-vanishing `lower order' invariants. In particular, investigating the problem of homotoping four 2-spheres to be disjoint in a simply-connected 4-manifold, they were led to a set of lattice points in $\mathbf{Z}^2,$ and asked if this set is an additive subgroup of $\mathbf{Z}^2.$  Projecting this set onto its first coordinate gives precisely the set of integers considered in Theorem~\ref{knsumproduct:theorem:subgroup}, and Theorem~\ref{knsumproduct:theorem:subgroup} states that this set is a subgroup of the additive group of integers.

\section{Sums of products of pairs}
\bl        \label{knsumprod:lemma:key}
Let $a,b,c,d,$ and $m$ be integers with $m \geq 1,$ and let 
\beq   \label{Nab}
N \equiv ab+cd \pmod{m}.
\eeq
If there exist integers $a'$ and $c'$ such that 
\[
a'\equiv a \pmod{m} 
\]
\[
 c'\equiv c\pmod{m}
\]
\[
\gcd(a',c')=m'
\]
and
\[
N \equiv a'b+c'd \pmod{mm'}
\]
then there exist integers $b'$ and $d'$ such that 
\beq  \label{knsumprod:b'}
b'\equiv b \pmod{m}
\eeq
\beq  \label{knsumprod:d'}
d'\equiv d \pmod{m}
\eeq
\beq  \label{knsumprod:N}
N = a'b'+c'd'.
\eeq
If $a'$ and $c'$ are positive and 
\beq  \label{N}
N  \geq a'b+c'd + m(a'-m')(c'-m')
\eeq
then there exist integers $b' \geq b$ and $d' \geq d$ that 
satisfy~\eqref{knsumprod:b'},~\eqref{knsumprod:d'}, 
and ~\eqref{knsumprod:N}.
\el

\begin{proof}
Congruence~\eqref{Nab} implies that there is an integer $\ell$ such that 
\[
N= a'b+c'd + \ell mm'.
\]
Since $m'= \gcd(a',c'),$ there exist integers $r$ and $s$ such that
\[
a'r + c' s = \ell m' = \frac{N-(a'b+c'd)}{m}.
\]
Defining $b'=b+mr$ and $d' = d+ms,$ we obtain
\[
N = a'b'+c'd'.
\]
A theorem of Sylvester~\cite[Theorem 1.17]{nath00aa}, which is a special case of the linear diophantine problem of Frobenius, implies that if $a'$ and $c'$ are positive integers with $(a',c')=m'$ and if 
\[
\ell \geq \left(\frac{a'}{m'} - 1 \right)\left(\frac{c'}{m'} - 1 \right)
\]
then there exist nonnegative integers $r$ and $s$ such that $a'r + c' s = \ell m' .$  It follows that if $N$ satisfies inequality~\eqref{N}, then $b' \geq b$ and $d' \geq d$.  This completes the proof.
\end{proof}

\bt   \label{knsumprod:theorem:main}
If $a,b,c,d,$ and $m$ are integers with $m \geq 1$ and $(a,b,c,d,m)=1,$  then 
\[
R_m(a)R_m(b)+R_m(c)R_m(d) = R_m(ab+cd).
\]
\et

\begin{proof}
Since we are only interested in the congruence classes of $a,b,c,d$ modulo $m$, we can assume without loss of generality that $a,b,c,d$ are positive.

Let $N \in R_m(ab+cd).$  There is an integer $k$ such that
\[
N = ab+cd +km.
\]
We define
\[
m' = \gcd(a,c,m).
\]
Then $1 \leq m' \leq m$.  Since  
\[
\gcd(b,d,m')= \gcd(a,b,c,d,m) = 1
\]
there are integers $x,y,z$ such that  
\beq  \label{A}
bx+dy+m'z= k = \frac{N-(ab+cd)}{m}. 
\eeq
Choose integers $x'$ and $y'$ such that 
\[
x' \equiv x \pmod{m'} \qquad\text{and}\qquad 0 \leq x' \leq m'-1
\]
and
\[
y' \equiv y \pmod{m'} \qquad\text{and}\qquad  bm \leq y' \leq bm + m'-1.
\]
There are integers $q_x$ and $q_y$ such that 
\beq  \label{B}
x=q_x m' + x' \qquad \text{and}  \qquad y = q_y m' + y'.
\eeq
It follows from~\eqref{A} and~\eqref{B} that
\begin{align*}
N & = (a+mx)b+(c+my)d+mm'z \\
& = (a+mx')b+(c+my')d+mm'(bq_x + dq_y + z).
\end{align*}
Let
\[
a_0 = a+mx' \quad\text{ and } \quad c_0= c+my'.
\]
Then $a_0 \geq 1$ since $a \geq 1, m \geq 1,$ and $x' \geq 0.$  We have 
\begin{align*}
a_0 & \equiv a \pmod{m}   \\
c_0 & \equiv c \pmod{m}  \\
\gcd(a_0,c_0,m) & = \gcd(a,c,m)=m'
\end{align*}
\[
N\equiv a_0b+c_0d\pmod{mm'}
\]
and
\beq     \label{ineq-0}
\left\{
\begin{array}{rcccl}
a & \leq & a_0 & < & a+m^2  \\
c+bm^2  & \leq & c_0& < & c+(b+1)m^2.
\end{array}
\right.
\eeq
Since $m'$ divides $\gcd(a_0,c_0),$ we have 
\[
\ord_p(m') \leq \ord_p(\gcd(a_0,c_0))
\]
for all prime numbers $p$.  

Let $\mathcal{P}$ be the set of prime numbers that divide $m'.$  The set $\mathcal{P}$ is finite because $m'\neq 0.$  Then $\mathcal{P} = \mathcal{P}_1 \cup \mathcal{P}_2,$ where
\[
\mathcal{P}_1 = \{ p \in \mathcal{P} : \ord_p(m') < \ord_p(\gcd(a_0,c_0)) \}
\]
\[
\mathcal{P}_2 = \{ p \in \mathcal{P} : \ord_p(m') = \ord_p(\gcd(a_0,c_0)) \}
\]
and $\mathcal{P}_1 \cap \mathcal{P}_2 = \emptyset.$
By the Chinese remainder theorem, there is an integer $u$ such that 
\begin{align*}
u \equiv 1 \pmod{p} & \qquad\text{for all $p \in \mathcal{P}_1$} \\
u \equiv 0 \pmod{p} & \qquad\text{for all $p \in \mathcal{P}_2$} 
\end{align*}
and
\[
0 \leq u < \prod_{p\in \mathcal{P}} p \leq m' \leq m.
\]
We define 
\begin{align*}
a_1 & = a_0+dmu \\
c_1& = c_0-bmu.
\end{align*}
Then
\[
(a_1,c_1,m)=(a_0,c_0,m)=m'
\]
and so 
\[
\ord_p(m') \leq  \ord_p(\gcd(a_1,c_1)).
\]
Since
\[
a_1b+c_1d = a_0b+c_0d
\]
we have
\[
N\equiv a_1 b+c_1 d\pmod{mm'}.
\]
Inequality~\eqref{ineq-0} implies that
\beq   \label{ineq-1}
\left\{
\begin{array}{rcl}
a \leq a_0  \leq & a_1 & \leq  a_0 + dm^2 \leq a + (d+1)m^2  \\
c \leq c_0-bm^2  \leq & c_1 & \leq  c_0 \leq c + (b+1)m^2
\end{array}
\right.
\eeq

Let $p \in \mathcal{P}_1$.  Then $p$ does not divide $u.$  Moreover, $m' = \gcd(a_0,c_0,m)$ and $\ord_p(m') < \ord_p(\gcd(a_0,c_0))$ implies that $\ord_p(m) = \ord_p(m')$.  Since $p$ divides $m'$ and $\gcd(b,d,m')=1,$ it follows that either $p$ does not divide $b$ or $p$ does not divide $d.$   If  $p$ does not divide $d,$   then
\[
\ord_p(dmu) =\ord_p(m) = \ord_p(m') < \ord_p(\gcd(a_0,c_0)) \leq \ord_p(a_0)
\]
and so
\[
\ord_p(a_1)=\ord_p(dmu) = \ord_p(m').
\]
Similarly, if $p$ does not divide $b,$  then 
\[
\ord_p(c_1)=\ord_p(bmu) = \ord_p(m').
\]
It follows that
\[
\ord_p(m') \geq \min\left(\ord_p(a_1),\ord_p(c_1)  \right)
= \ord_p(\gcd(a_1,c_1))
\]
and so
\[
\ord_p(m') = \ord_p(\gcd(a_1,c_1))\qquad\text{for all $p \in \mathcal{P}_1.$}
\]

Let $p \in \mathcal{P}_2.$  Then $p$ does divide $u.$  Since $\ord_p(m') = \ord_p(\gcd(a_0,c_0)),$ it follows that either
\[
\ord_p(m') = \ord_p(a_0) \leq \ord_p(c_0)
\]
or 
\[
\ord_p(m') = \ord_p(c_0) \leq \ord_p(a_0).
\]
In the first case,
\[
\ord_p(a_0) = \ord_p(m') \leq \ord_p(m) < \ord_p(dmu)
\]
and so $\ord_p(a_1) = \ord_p(a_0) = \ord_p(m') \leq \ord_p(c_1).$
In the second case, $\ord_p(c_1) = \ord_p(m') \leq \ord_p(a_1).$
Therefore,
\[
\ord_p(m') = \ord_p(\gcd(a_1,c_1))\qquad\text{for all $p \in \mathcal{P}_2.$}
\]
It follows that
\[
\ord_p(m') = \ord_p(\gcd(a_1,c_1))\qquad\text{for all $p \in \mathcal{P}.$}
\]

Let $\mathcal{P}_3$ be the set of prime numbers that divide $a_1$ but do not divide $m$.   Then $\mathcal{P}_3$ is finite since $a_1\neq 0.$  
By the Chinese remainder theorem, since $\gcd(mm',p)=1$, there is an integer $v$ such that 
\[
0 \leq v \leq \prod_{p\in \mathcal{P}_3} p \leq a_1
\]
and
\[
c_1+mm'v \equiv 1 \pmod{p}\qquad\text{for all $p \in \mathcal{P}_3$.}
\]
Let $a'=a_1$ and $c'=c_1+mm'v.$  
Inequality~\eqref{ineq-1} implies that
\beq   \label{ineq-2}
\left\{
\begin{array}{rcl}
a \leq & a' & \leq  a + (d+1)m^2  \\
c \leq c_1 \leq & c' & \leq c_1+ a_1m^2 \leq c+(a+b+1)m^2 +(d+1)m^4
\end{array}
\right.
\eeq
If a prime $p$ divides $a'$ but does not divide $m$, then $p$ does not divide $c'.$   Thus, if a prime divides both $a'$ and $c',$ then it must divide $m$ and so it divides $m'.$  Since $\gcd(a',c') = \gcd(a_1,c_1),$ it follows that if $p$ divides both $a'$ and $c'$, then $\ord_p(\gcd(a',c')) = \ord_p(m')$ and so
\[
\gcd(a',c') = m'.
\]
We also have 
\[
a'\equiv a \pmod{m}
\]
\[
c'\equiv c \pmod{m}
\]
\[
N \equiv a'b+c'd\pmod{mm'}.
\]
Lemma~\ref{knsumprod:lemma:key} implies that $N \in R_m(a)R_m(b)+R_m(c)R_m(d).$
This completes the proof.
\end{proof}

\bt     \label{knsumprod:theorem:main-delta}
Let $a,b,c,d,$ and $m$ be integers with $m \geq 1$ and $(a,b,c,d,m)=\delta.$  Then 
\[
R_m(a)R_m(b)+R_m(c)R_m(d) = R_{\delta m}(ab+cd).
\]
\et

\begin{proof}
For all integers $A,M,$ and $\delta$ with $M \geq 1$ and $\delta \geq 1$ we have
\[
\delta \ast R_M(A) = R_{\delta M}(\delta A).
\]
Let $A = a/\delta, B= b/\delta, C=c/\delta, D=d/\delta,$ and $M = m/\delta.$  Since $\gcd(A,B,C,D,M)=1$ and $M \geq 1,$  Theorem~\ref{knsumprod:theorem:main} implies that 
\[
R_M(A)R_M(B) + R_M(C)R_M(D) = R_M(AB+CD).
\]
Multiplying by $\delta^2,$ we obtain
\begin{align*}
R_m(a)R_m(b)+R_m(c)R_m(d)
& = \delta^2  \ast \left( R_M(A)R_M(B) + R_M(C)R_M(D) \right) \\
& = \delta^2 R_M(AB+CD) \\
& =  R_{\delta^2  \ast M}( (\delta  A)( \delta  B) + (\delta  C) ( \delta  D) )\\
& = R_{\delta m}(ab+cd).
\end{align*}
This completes the proof.
\end{proof}

\bt     \label{knsumproduct:theorem:subgroup}
Let $a,b,c,d,$ and $m$ be integers with $m \geq 1$ and $(a,b,c,d,m)=\delta.$  Then 
\[
\{aw+bx+cy+dz+m(wx+yz) : w,x,y,z \in \Z\} = \delta \Z.
\]
\et

\begin{proof}
This is simply an unraveling of Theorem~\ref{knsumprod:theorem:main-delta}.
\end{proof}

\bt   \label{knsumprod:theorem:mainAP}
If $a,b,c,d,$ and $m$ are positive integers with $(a,b,c,d,m)=1,$  then 
\[
P_m(a)P_m(b)+P_m(c)P_m(d) \sim P_m(ab+cd).
\]
\et

\begin{proof}
Let $N \in P_m(ab+cd)$.  We must prove that if $N$ is sufficiently large, then $N \in P_m(a)P_m(b)+P_m(c)P_m(d).$  In the proof of Theorem~\ref{knsumprod:theorem:main} we constructed positive integers $a'$ and $c'$ satisfying inequalities~\eqref{ineq-2}
\[
a \leq a' \leq  a + (d+1)m^2
\]
\[
c \leq c' \leq  c+(a+b+1)m^2 +(d+1)m^4
\]
and the hypotheses of Lemma~\ref{knsumprod:lemma:key}:
\[
a'\equiv a \pmod{m} 
\]
\[
 c'\equiv c\pmod{m}
\]
\[
\gcd(a',c')=m'
\]
and
\[
N \equiv a'b+c'd \pmod{mm'}
\]
Let 
\begin{align*}
N_0 = & (a + (d+1)m^2)b + (c+(a+b+1)m^2 +(d+1)m^4)d \\
& + m  (a + (d+1)m^2) (c+(a+b+1)m^2 +(d+1)m^4).
\end{align*}
Then 
\[
N_0 \geq a'b+c'd+ma'c' > a'b+c'd+m(a'-m')(c'-m').
\]
It follows from Lemma~\ref{knsumprod:lemma:key} that if $N \geq N_0,$ then there exist integers $b' \in P_m(b)$ and $d' \in P_m(d)$ such that $N=a'b'+c'd'.$  This completes the proof.
\end{proof}

\section{Iterated sums and products}
Let $h \geq 2$ and let $k_1,k_2,\ldots,k_h$ be positive integers such that $k_1\leq k_2 \leq \cdots \leq k_h.$  Let $a_{i,j}\in \Z$ for $ i=1,\ldots, h$ and $j=1,\ldots,k_i$, and let $m \geq 1.$   Then
\[
\sum_{i=1}^h \prod_{j=1}^{k_i} R_m(a_{i,j}) \subseteq R_m\left( \sum_{i=1}^h \prod_{j=1}^{k_i} a_{i,j} \right).
\]
We would like to know when the inclusion is an equality, that is, when we have
\[
\sum_{i=1}^h \prod_{j=1}^{k_i} R_m(a_{i,j}) = R_m\left( \sum_{i=1}^h \prod_{j=1}^{k_i} a_{i,j} \right).
\]

\bl   \label{knsumprod:lemma:k1=1}
Let $a_0,a_1,a_2, \ldots,a_k,$ and $m$ be integers with $m\geq 1.$  Then
\[
R_m(a_0) + R_m(a_1) R_m(a_2) \cdots  R_m(a_k) = R_m(a_0+a_1a_2\cdots a_k).
\]
\el
 
\begin{proof}
If 
\[
N \in R_m(a_0) + R_m(a_1) R_m(a_2) \cdots  R_m(a_k)
\]
then there are integers $q_0,q_1,q_2,\ldots, q_k$ such that
\[
N = (a_0+q_0m) + (a_1 + q_1 m)(a_2 + q_2 m) \cdots (a_k + q_k m)
\]
and so
\[
N \equiv a_0+a_1a_2\cdots a_k \pmod{m}.
\]
It follows that
\[
R_m(a_0) + R_m(a_1) R_m(a_2) \cdots  R_m(a_k) \subseteq R_m(a_0+a_1a_2\cdots a_k).
\]
Conversely, if 
\[
N \in R_m(a_0+a_1a_2\cdots a_k)
\]
then there is an integer $q$ such that 
\begin{align*}
N 
& =  a_0+a_1a_2\cdots a_k + qm \\
& = (a_0+qm) + a_1a_2\cdots a_k\\
& \in R_m(a_0) + R_m(a_1) R_m(a_2) \cdots  R_m(a_k)
\end{align*}
and so 
\[
R_m(a_0+a_1a_2\cdots a_k) \subseteq  R_m(a_0) + R_m(a_1) R_m(a_2) \cdots  R_m(a_k).
\]
This completes the proof.
\end{proof}

\bt
Let $h \geq 2$ and let $k_1,k_2,\ldots,k_h$ be positive integers such that $k_1\leq k_2 \leq \cdots \leq k_h.$  Let $a_{i,j}\in \Z$ for $ i=1,\ldots, h$ and $j=1,\ldots,k_i$, and let $m \geq 1.$   If $k_1=1$ or if $k_1=k_2=2$ and $\gcd(a_{1,1},a_{1,2},a_{2,1},a_{2,2},m)=1,$ then 
\[
\sum_{i=1}^h \prod_{j=1}^{k_i} R_m(a_{i,j}) = R_m\left( \sum_{i=1}^h \prod_{j=1}^{k_i} a_{i,j} \right).
\]
\et

\begin{proof}
The case $h=2$ follows immediately from Lemma~\ref{knsumprod:lemma:k1=1} and Theorem~\ref{knsumprod:theorem:main}, and the result for all $h\geq 2$ follows by induction.
\end{proof}

\def\cprime{$'$} \def\cprime{$'$} \def\cprime{$'$}
\providecommand{\bysame}{\leavevmode\hbox to3em{\hrulefill}\thinspace}
\providecommand{\MR}{\relax\ifhmode\unskip\space\fi MR }
\providecommand{\MRhref}[2]{%
  \href{http://www.ams.org/mathscinet-getitem?mr=#1}{#2}
}
\providecommand{\href}[2]{#2}

\end{document}